\documentclass{amsart}
\usepackage{epsf}
\usepackage{amscd}             

\def\Mapright#1#2{\,\smash{\mathop{\smash{\hbox to #1pt{\rightarrowfill}}
\vphantom\rightarrow}\limits^{#2}\,}}
\def\co{\hbox{$\,:\,$}} 
\def\larrow{\mathop{\longrightarrow}\limits}
\newcommand{\Jscr}{\mathcal J}
\newcommand{\Hscr}{\mathcal H}
\newcommand{\vol}{\operatorname{vol}}
\newcommand{\Ker}{\operatorname{Ker}}
\newcommand{\Hom}{\operatorname{Hom}}
\renewcommand{\Im}{\operatorname{Im}}

\newcommand{\cs}{\operatorname{cs}}
\newcommand{\Isom}{\operatorname{Isom}}
\newcommand{\PSL}{\operatorname{PSL}}
\newcommand\calE{{\mathcal E}}
\newcommand\calR{{\mathcal R}}
\newcommand\eeprebloch{{\calE\prebloch}}
\newcommand\eebloch{{\calE\bloch}}
\newcommand\prebloch{{\mathcal P}}
\newcommand\bloch{{\mathcal B}}
\newcommand\FT{\operatorname{FT}}
\newcommand\liftFT{\widehat\FT}
\newcommand\eprebloch{{\widehat\prebloch}}
\newcommand\ebloch{{\widehat\bloch}}
\newcommand\C{\mathbb C}

\newcommand\Q{\mathbb Q}
\newcommand\CP{\C P}
\newcommand\Z{\mathbb Z}
\newcommand\bfw{{\bf w}}
\renewcommand\H{\mathbb H}
\newcommand\discrete{^\delta}
\newcommand\abc{(\C-\{0,1\})^{~ab}}
\newcommand\themap{\lambda}
\newtheorem{theorem}{Theorem}[section]
\newtheorem{lemma}[theorem]{Lemma}
\newtheorem{proposition}[theorem]{Proposition}
\newtheorem{corollary}[theorem]{Corollary}
\theoremstyle{definition}
\newtheorem{definition}[theorem]{Definition}
\newcommand\Cover{\widehat\C}
\begin{document}
\title[Extended Bloch group and the Chern-Simons class]{Extended Bloch group 
and the Chern-Simons class\\(Incomplete Working version)}
\author{Walter D. Neumann}
\address{Department of Mathematics\\The University of
Melbourne\\Carlton, Vic 3052\\Australia}
\email{neumann@maths.mu.oz.au}

\begin{abstract}
We define an extended Bloch group and show it is isomorphic 
to $H_3(\PSL(2,\C)^\delta;\Z)$. Using the Rogers dilogarithm function
this leads to an exact simplicial formula for the universal
Cheeger-Simons class on this homology group.  It also leads to an
independent proof of the analytic relationship between volume and
Chern-Simons invariant of hyperbolic manifolds conjectured in
\cite{neumann-zagier} and proved in \cite{yoshida}, as well as an
effective formula for the Chern-Simons invariant of a hyperbolic
manifold.
\end{abstract}
\maketitle

\section{Introduction} 
There are several variations of the definition of the Bloch group in
the literature; by \cite{dupont-sah} they differ at most by torsion
and they agree with each other for algebraically closed fields.  In
this paper we shall use the following.

\begin{definition}\label{def-bloch}
Let $k$ be a field.  The \emph{pre-Bloch group $\prebloch(k)$} is the
quotient of the free $\Z$-module $\Z (k-\{0,1\})$ by all instances of
the following relation:
\begin{equation}
[x]-[y]+[\frac yx]-[\frac{1-x^{-1}}{1-y^{-1}}]+[\frac{1-x}{1-y}]=0,
\label{5term}
\end{equation}
This relation is usually called the \emph{five term relation}.  The
\emph{Bloch group $\bloch(k)$} is the kernel of the map
\begin{equation*}
\prebloch(k)\to k^*\wedge_\Z k^*,\quad [z]\mapsto 2\bigl(z\wedge(1-z)\bigr).
\end{equation*}
\end{definition}

(In \cite{neumann-yang} the additional relations
\begin{equation*}
[x]=[1-\frac 1x]=[\frac 1{1-x}]=-[\frac1x]=-[\frac{x-1}x]=-[1-x]
\label{invsim}
\end{equation*}
were used.  These follow from the five term relation when $k$ is
algebraicly closed, as shown by
Dupont and Sah \cite{dupont-sah}.  Dupont and Sah use a different five
term relation but it is conjugate to the one used here by
$z\mapsto\frac1z$.)

There is an exact sequence due to Bloch and Wigner: 
\begin{equation*}
0\to \Q/\Z\to H_3(\PSL(2,\C)\discrete;\Z)\to \bloch(\C)\to 0.
\end{equation*}
The superscript $\delta$ means ``with discrete topology.''
\renewcommand\discrete{} We will omit it from now on.

$\bloch(\C)$ is known to be uniquely divisible, so it has canonically
the structure of a $\Q$-vector space (Suslin \cite{suslin}). It's
$\Q$-dimension is infinite and conjectured to be countable (the
``Rigidity Conjecture,'' equivalent to the conjecture that
$\bloch(\C)=\bloch(\overline\Q)$, where $\overline\Q$ is the field of
algebraic numbers).  In particular, the $\Q/\Z$ in the Bloch-Wigner
exact sequence is precisely the torsion of $H_3(\PSL(2,\C);\Z)$, so
any finite torsion subgroup is cyclic.

\medskip 
In the present paper we define an \emph{extended Bloch group}
$\eebloch(\C)$ by replacing $\C-\{0,1\}$ in the definition of
$\bloch(\C)$ by its universal abelian cover $\abc$ and
appropriately lifting the five term relation (\ref{5term}). Our main
results are that we can lift the Bloch-Wigner map
$H_3(\PSL(2,\C)\discrete;\Z)\to \bloch(\C)$ to an isomorphism
\begin{equation*}
\themap\colon H_3(\PSL(2,\C)\discrete;\Z)\to \eebloch(\C)
\end{equation*}
Moreover, the ``Roger's dilogarithm function'' (see below) gives a
natural map
\begin{equation*}R\colon\eebloch(\C)\to\C/2\pi^2\Z.\end{equation*} We
show that the composition \begin{equation*}R\circ\themap\colon
H_3(\PSL(2,\C)\discrete;\Z)\to\C/2\pi^2\Z\end{equation*} is the
Cheeger-Simons class (cf \cite{cheeger-simons}), so it can also be
described as $i(\vol + i \cs)$, where $cs$ is the universal
Chern-Simons class.  It has been a longstanding problem to provide
such a computation of the Chern-Simons class. Dupont in \cite{dupont1}
gave an answer modulo $\pi^2\Q$ and our computation is a natural lift
of his.

Another consequence of our result is that any complete hyperbolic
3-manifold $M$ of finite volume has a natural ``fundamental class'' in
$H_3(\PSL(2,\C)\discrete;\Z)/C_2$, where $C_2$ is the (unique) order 2
subgroup.  For compact $M$ the existence of this class, even without
the $C_2$ ambiguity, is easy and well known: $M=\H^3/\Gamma$ is a
$K(\Gamma,1)$-space, so the inclusion $\Gamma\to \Isom^+(\H^3) =
\PSL(2,\C)$ induces $H_3(M)=H_3(\Gamma)\to H_3(\PSL(2,\C)$, and the
class in question is the image of the fundamental class $[M]\in
H_3(M)$.  For non-compact $M$ the existence of such a class is
somewhat surprising, although it was already strongly suggested by
earlier results.

We can describe this fundamental class nicely in terms of an ideal
triangulation of $M$.  However, this ideal triangulation has to be a
``true'' ideal triangulation rather than the less restrictive ``degree
1'' ideal triangulations used in \cite{neumann-yang}.  The ideal
triangulations resulting from Dehn filling that are used by the
programs Snappea and Snap \cite{snap} are not true.
Nevertheless, we can describe the fundamental class in terms of these
``Dehn filling triangulations.'' This leads also to an exact
simplicial formula for the Chern-Simons invariant of a hyperbolic
3-manifold, refining the formula of \cite{neumann}.

We work initially with a different version of the extended Bloch
group, based on a disconnected $\Z\times\Z$ cover of $\C-\{0,1\}$.
This group, which we call $\ebloch(\C)$ is a quotient of
$\eebloch(\C)$ by a subgroup of order $2$. 

\medskip\noindent{\bf Acknowledgements.} The definition of the extended Bloch
group was suggested by an idea of Jun Yang, to whom I am grateful also
for many useful conversations.  In particular, he informs me that this
work can be interpreted as giving a motivic complex for $K_3(\C)$.
The main results of this paper were announced in
\cite{neumann-hilbert}. This research is supported by the Australian
Research Council.

\section{The preliminary version of extended Bloch group}

We shall need a $\Z\times\Z$ cover $\Cover$ of $\C-\{0,1\}$ which can
be constructed as follows.  Let $P$ be $\C-\{0,1\}$ split along the
rays $(-\infty,0)$ and $(1,\infty)$.  Thus each real number $r$
outside the interval $[0,1]$ occurs twice in $P$, once in the upper
half plane of $\C$ and once in the lower half plane of $\C$.  We
denote these two occurences of $r$ by $r+0i$ and $r-0i$.  We construct
$\Cover $ as an identification space from $P\times\Z\times\Z$ by
identifying
\begin{equation*}
\begin{aligned}
(x+0i, p,q)&\sim (x-0i,p+2,q)\quad\hbox{for each }x\in(-\infty,0)\\
(x+0i, p,q)&\sim (x-0i,p,q+2)\quad\hbox{for each }x\in(1,\infty).\\
\end{aligned}
\end{equation*}
We will denote the equivalence class of $(z,p,q)$ by $(z;p,q)$.
$\Cover $ has four components:
\begin{equation*}
\Cover =X_{00}\cup X_{01}\cup X_{10}\cup X_{11}
\end{equation*}
where $X_{\epsilon_0\epsilon_1}$ is the set of $(z;p,q)\in \Cover $ with
$p\equiv\epsilon_0$ and $q\equiv\epsilon_1$ (mod $2$).

We may think of $X_{00}$ as the riemann surface for the function
$\C-\{0,1\}\to\C^2$ defined by $z\mapsto \bigl(\log z, -\log (1-z)\bigr)$.
Taking the branch $(\log z + 2p\pi i, -\log (1-z) + 2q\pi i)$ of this
function on the portion $P\times\{(2p,2q)\}$ of $X_{00}$ for each
$p,q\in\Z$ defines an analytic function from $X_{00}$ to $\C^2$.  In
the same way, we may think of $\Cover $ as the riemann surface for the
collection of all branches of the functions $(\log z + p\pi i, -\log
(1-z) + q\pi i)$ on $\C-\{0,1\}$.

\medbreak
Consider the set
\begin{equation*}
\FT:=\biggl\{\biggl(x,y,\frac yx,\frac{1-x^{-1}}{1-y^{-1}},
\frac{1-x}{1-y}\biggr):x\ne y,
x,y\in\C-\{0,1\}\biggr\}\subset(\C-\{0,1\})^5
\end{equation*}
of 5-tuples involved in the five term relation (\ref{5term}).  
An elementary computation shows:

\begin{lemma}
The subset $\FT^+$ of $(x_0,\dots,x_4)\in \FT$ with each $x_i$ 
in the upper half plane of $\C$ is the set of elements of $\FT$ for
which $y$ is in the upper half plane of $\C$ and $x$ is inside the
triangle with vertices $0,1,y$. Thus $\FT^+$ is connected (even
contractible).  \qed\end{lemma}

\begin{definition}\label{liftFT}
Let $V\subset(\Z\times\Z)^5$ be the subspace
\begin{multline*}
V:=\{\bigl((p_0,q_0),(p_1,q_1),(p_1-p_0,q_2),
(p_1-p_0+q_1-q_0,q_2-q_1),\\ (q_1-q_0,q_2-q_1-p_0)\bigr):p_0,p_1,q_0,q_1,q_2\in\Z\}.
\end{multline*}
Let $\liftFT_0$ denote the unique component of the inverse image of $\FT$ in
$\Cover^5$ which includes the points $\bigl((x_0;0,0),\ldots,(x_4;0,0)\bigr)$
with $(x_0,\ldots,x_4)\in \FT^+$, and define
\begin{equation*}
\liftFT:=\liftFT_0+V=\{
{\mathbf x}+{\mathbf v}:{\mathbf x}\in
\liftFT_0\text{ and }{\mathbf v}\in V\},
\end{equation*}
where we are using addition to denote the action of $(\Z\times\Z)^5$
by covering transformations on $\Cover^5$. (Although we do not need
it, one can show that the action of $2V$ takes $\liftFT_0$ to itself,
so $\liftFT$ has $2^5$ components, determined by the parities of
$p_0,p_1,q_0,q_1,q_2$.)

Define $\eprebloch(\C)$ as the free
$\Z$-module on $\Cover$ factored by all instances of the relations:
\begin{equation}
\sum_{i=0}^4(-1)^i(x_i;p_i,q_i)=0\label{five term}\quad
\text
{with $\bigl((x_0;p_0,q_0),\dots,(x_4;p_4,q_4)\bigr)\in\liftFT$}
\end{equation}
and
\begin{equation}
(x;p,q)+(x;p',q')=(x;p,q')+(x;p',q)\label{transfer}
\quad\text{with $p,q,p',q'\in\Z$}
\end{equation}
We shall denote the class of $(z;p,q)$ in
$\eprebloch(\C)$ by $[z,p,q]$
\end{definition}

We call relation (\ref{five term}) the \emph{lifted five term relation}.
We shall see that its precise form arises naturally in several
contexts. In particular, we give it a geometric interpretation in
Sect.~\ref{simplex parameters}.  

We call relation (\ref{transfer}) the \emph{transfer relation}.  It is
almost a consequence of the lifted five term relation, since we shall
see that the effect of omitting it would be to replace
$\eprebloch(\C)$ by $\eprebloch(\C)\otimes\Z/2$, with $\Z/2$ generated
by an element $\kappa:=[x,1,1]+[x,0,0]-[x,1,0]-[x,0,1]$ which is
independent of $x$.

\begin{lemma}
There is a well-defined homomorphism
\begin{equation*}\nu\colon\eprebloch(\C)\to \C\wedge_\Z\C
\end{equation*} defined on generators by
$[z,p,q]\mapsto (\log z + p\pi i)\wedge (-\log(1-z) +q\pi i)$.  
\end{lemma}
\begin{proof}
We must verify that $\nu$ vanishes on the relations that define
$\eprebloch(\C)$. This is trivial for the transfer relation
(\ref{transfer}).
We shall show that the lifted five term
relation is the most general lift of the five term relation
(\ref{5term}) for which $\nu$ vanishes.  If one applies $\nu$ to an
element $\sum_{i=0}^4(-1)^i[x_i,p_i,q_i]$ with
$(x_0,\ldots,x_4)=(x,y,\ldots)\in \FT^+$ one obtains after simplification:
\begin{multline*}
\bigl((q_0-p_2-q_2+p_3+q_3)\log x
+(p_0-q_3+q_4)\log(1-x)+(-q_1+q_2-q_3)\log y +{}\\ {}+
(-p_1+p_3+q_3-p_4-q_4)\log(1-y)+(p_2-p_3+p_4)\log(x-y)\bigr)\wedge \pi i.
\end{multline*}
An elementary linear algebra computation shows that this vanishes
identically if and only if $p_2=p_1-p_0$, $p_3=p_1-p_0+q_1-q_0$,
$q_3=q_2-q_1$, $p_4=q_1-q_0$, and $q_4 = q_2-q_1-p_0$, as in the
lifted five term relation.  The vanishing of $\nu$ for the general
lifted five term relation now follows by analytic continuation.
\end{proof}
\begin{definition}
Define $\ebloch(\C)$ as the kernel of
$\nu\colon\eprebloch(\C)\to\C\wedge\C$.
\end{definition}
 
Define
\begin{equation*}
R(z;p,q)={\calR}(z)+\frac{\pi i}2(p\log(1-z)+q\log z) -\frac{\pi^2}6
\end{equation*} 
where $\calR$ is the Rogers dilogarithm function
\begin{equation*}
{\calR}(z)=\frac12\log(z)\log(1-z)-\int_0^z\frac{\log(1-t)}tdt.\end{equation*}
Then

\begin{proposition}\label{R}
$R$ gives a well defined map $R\colon \Cover\to \C/\pi^2\Z$.  The relations
which define $\eprebloch(\C)$ are functional equations for $R$ modulo
$\pi^2$ (the lifted five term relation is in fact the most general
lift of the five term relation (\ref{5term}) with this property).
Thus $R$ also gives a homomorphism
$R\colon\eprebloch(\C)\to\C/\pi^2\Z$.
\end{proposition}

\begin{proof}
If one follows a closed path from $z$
that goes anti-clockwise around the origin it is easily verified that
$R(z;p,q)$ is replaced by $R(z;p,q)+\pi i \log(1-z)
-q\pi^2=R[z,p+2,q]-q\pi^2$.  Similarly, following a closed path
clockwise around $1$ replaces $R(z;p,q)$ by $R(z;p,q+2)+p\pi^2$.  Thus
$R$ modulo $\pi^2$ is well defined on $\Cover$ (in fact $R$ itself is well
defined on a $\Z$ cover of $\Cover$ which is a nilpotent cover of
$\C-\{0,1\}$). 

It is well known that
$L(z):=\calR(z)-\frac{\pi^2}6$ satisfies the functional equation
\begin{equation*}
L(x)-L(y)+L\bigl(\frac yx\bigr)-L\bigl(\frac{1-x^{-1}}{1-y^{-1}}\bigr)+
L\bigl(\frac{1-x}{1-y}\bigr)=0\end{equation*} 
for
$0<y<x<1$.  Since the 5-tuples involved in this equation are on the
boundary of $\FT^+$, the functional equation
\begin{equation*}\sum(-1)^iR(x_i;0,0)=0\end{equation*}
is valid by analytic continuation on the whole of $\FT^+$.  Now
\begin{equation*}\sum(-1)^iR(x_i;p_i,q_i)\end{equation*}
differs from this by 
\begin{equation*}\frac{\pi i}2\sum(-1)^i(p_i\log(1-x_i)+q_i\log x_i)
\end{equation*} 
and it is an elementary calculation to verify that this vanishes
identically for $(x_0,\ldots,x_4)\in \FT^+$ if and only if the $p_i$
and $q_i$ are as in the lifted five term relation.  Thus the lifted
five-term relation gives a functional equation for $R$ when
$(x_0,\ldots,x_4)\in \FT^+$.  By analytic continuation, it is a
functional equation for $R$ mod $\pi^2$ in general.  The transfer
relation is trivially a functional equation for $R$.
\end{proof}

The first version of our main result is
\begin{theorem}
There exists an epimorphism $\themap\colon
H_3(\PSL(2,\C)\discrete;\Z)\to\ebloch(\C)$ with kernel of
order 2 such that the composition $\themap\circ R\colon
H_3(\PSL(2,\C)\discrete;\Z)\to\C/\pi^2\Z$ is the characteristic class
given by $i(\vol + i \cs)$.
\end{theorem}

We shall later modify the definition of $\ebloch(\C)$ to eliminate the
$\Z/2$ kernel.  To describe the map $\themap$ we must
give a geometric interpretation of $\Cover$.
 
\section{Parameters for ideal hyperbolic simplices}\label{simplex
parameters}

In this section we shall interpret $\Cover$ as a space of parameters for
what we call ``combinatorial flattenings'' of ideal hyperbolic
simplices.  We need this to define the above map $\themap$.  It also
gives a geometric interpretation of the lifted five term relation.

We shall denote the standard compactification of $\H^3$ by $\overline
\H^3 = \H^3\cup\CP^1$. An ideal simplex $\Delta$ with vertices
$z_1,z_2,z_3,z_4\in\CP^1$ is determined up to congruence by the cross
ratio
\begin{equation*}
z=[z_1\co z_2\co z_3\co z_4]=\frac{(z_3-z_2)(z_4-z_1)}{(z_3-z_1)(z_4-z_2)}.
\end{equation*}
Permuting the vertices by an even (i.e.,
orientation preserving) permutation replaces $z$ by one of
\begin{equation*}
z,\quad z'=\frac 1{1-z}, \quad\text{or}\quad z''=1-\frac 1z.
\end{equation*}
The parameter $z$ lies in the upper half plane of $\C$ if the
orientation induced by the given ordering of the vertices agrees with
the orientation of $\H^3$.  But we allow simplices whose vertex
ordering does not agree with their orientation. We also allow
degenerate ideal simplices whose vertices lie in one plane, so the
parameter $z$ is real.  However, we always require that the vertices
are distinct.  Thus the parameter $z$ of the simplex lies in
$\C-\{0,1\}$ and every such $z$ corresponds to an ideal simplex.

There is another way of describing the cross-ratio parameter $z=
[z_1\co z_2\co z_3\co z_4]$ of a simplex. The group of orientation preserving
isometries of $\H^3$ fixing the points $z_1$ and $z_2$ is
isomorphic to $\C^*$ and the element of this $\C^*$ that takes $z_4$
to $z_3$ is $z$. Thus the cross-ratio parameter $z$ is associated with
the edge $z_1z_2$ of the simplex.  The parameter associated in this
way with the other two edges $z_1z_4$ and $z_1z_3$ out of $z_1$ are
$z'$ and $z''$ respectively, while the edges $z_3z_4$, $z_2z_3$, and
$z_2z_4$ have the same parameters $z$, $z'$, and $z''$ as their
opposite edges. See fig.\ 1.
\begin{figure}[htbp]
\centerline{
\epsfxsize.35\hsize\epsffile{ebfig1b.eps}}
\caption{}
\end{figure}

Note that $zz'z''=-1$, so the sum
\begin{equation*}\log z + \log z' + \log z''\end{equation*} 
is an odd multiple of $\pi i$, depending on the branches of $\log$
used.  In fact, if we use standard branch of log then this sum is $\pi
i$ or $-\pi i$ depending on whether $z$ is in the upper or lower half
plane.
\begin{definition}\label{flattening}
We shall call any triple of the form
\begin{equation*}
\bfw=(w_0,w_1,w_2)=(\log z +p\pi i, \log z' + q\pi i, \log z'' + r\pi i)
\end{equation*} 
with
\begin{equation*}p,q,r\in \Z\quad\text{and}\quad
w_0+w_1+w_2=0
\end{equation*}
a \emph{combinatorial flattening} for our simplex. 

Each edge $E$ of $\Delta$ is assigned one of the
components $w_i$ of $\bfw$, with opposite edges being assigned the
same component. We call
$w_i$ the \emph{log-parameter} for the edge $E$ and denote it
$l_E(\Delta,\bfw)$. 
\end{definition}
This combinatorial flattening can be written 
\begin{equation*}\ell(z;p,q):=(\log z + p\pi i, -\log(1-z) + q\pi i,
\log(1-z)-\log z - (p+q)\pi i),
\end{equation*} 
and $\ell$ is then a map
of $\Cover$ to the set of combinatorial flattenings of simplices.

\begin{lemma}\label{Cover is flattenings} 
This map $\ell$ is a bijection, so $\Cover$ may be identified with
the set of all combinatorial flattenings of ideal tetrahedra.
\end{lemma}

\begin{proof}
We must show that we can recover $(z;p,q)$ from $(w_0,w_1,w_2)=
\ell(z;p,q)$.  It clearly suffices to recover $z$. But $z=\pm e^{w_0}$
and $1-z=\pm e^{-w_1}$, and the knowledge of both $z$ and $1-z$ up to
sign determines $z$.
\end{proof}

We can give a geometric interpretation of the choice of parameters in
the five term relation (\ref{five term}).  If $z_0,\ldots,z_4$ are
five distinct points of $\partial\overline\H^3$, then each choice of
four of five points $z_0,\dots,z_4$ gives an ideal simplex. We denote
the simplex which omits vertex $z_i$ by $\Delta_i$. The cross ratio
parameters $x_i=[z_0\co \ldots\co \hat{z_i}\co \ldots\co z_4]$ of
these simplices can be expressed in terms of $x:=x_0$ and $y:=x_1$ as
\begin{align*}
x_0=[z_1\co z_2\co z_3\co z_4]&=:x\\
x_1=[z_0\co z_2\co z_3\co z_4]&=:y\\
x_2=[z_0\co z_1\co z_3\co z_4]&=\frac yx\\
x_3=[z_0\co z_1\co z_2\co z_4]&=\frac {1-x^{-1}}{1-y^{-1}}\\
x_4=[z_0\co z_1\co z_2\co z_3]&=\frac {1-x}{1-y}
\end{align*}
The lifted
five term relation has the form \def\term#1{(x_#1;p_#1,q_#1)}
\begin{equation}\label{general 5-term}
\sum_{i=0}^4(-1)^i \term i=0
\end{equation}
with certain relations on the $p_i$ and $q_i$.  We will give a
geometric interpretation of
these relations.

Using the map of Lemma
\ref{Cover is flattenings}, each summand in this relation
(\ref{general 5-term}) represents a choice $\ell\term i$ of combinatorial
flattening for one of the five ideal simplices. For each edge $E$
connecting two of the points $z_i$ we get a corresponding linear
combination
\begin{equation}\label{edge sums}
\sum_{i=0}^4(-1)^il_E(\Delta_i,\ell\term i)
\end{equation}
of log-parameters (Definition \ref{flattening}), where we put
$l_E(\Delta_i,\ell\term i)=0$ if the line $E$ is not an edge of
$\Delta_i$.  This linear combination has just three non-zero terms
corresponding to the three simplices that meet at the edge $E$. One
easily checks that the real part is zero and the imaginary part can be
interpreted (with care about orientations) as the sum of the
``adjusted angles'' of the three flattened simplices meeting at $E$.

\begin{definition}
We say that the $\term i$ satisfy the \emph{flattening
condition} if each of the above linear combinations (\ref{edge sums})
of log-parameters is equal to zero. That is, the adjusted angle sum of
the three simplices meeting at each edge is zero.
\end{definition}

\begin{lemma}\label{geometric five term}
Relation (\ref{general 5-term}) is an instance of the lifted five term
relation (\ref{five term}) if and only if the $\term i$ satisfy
the flattening condition.
\end{lemma}
\begin{proof}
We first consider the case that $(x_0,\ldots,x_4)\in \FT^+$.  Recall
this means that each $x_i$ is in $\H$.  Geometrically, this implies
that each of the above five tetrahedra is positively oriented by the
ordering of its vertices.  This implies the configuration of fig.\ 2
\begin{figure}[hbtp]
\centerline{\epsfxsize.5\hsize\epsffile{ebfig2.eps}}
\caption{}
\end{figure}
with $z_1$ and $z_3$ on opposite sides of the plane of the triangle
$z_0z_2z_4$ and the line from $z_1$ to $z_3$ passing through the
interior of this triangle.  Denote the combinatorial flattening of the
$i^{th}$ simplex by $\ell(x_i;p_i,q_i)$.  If we consider the
log-parameters at the edge $z_3z_4$ for example, they are $\log x +
p_0\pi i$, $\log y + p_1 \pi i$, and $\log(y/x) + p_2\pi i$ and the
condition is that $(\log x + p_0\pi i)-(\log y + p_1
\pi i) +(\log(y/x) +p_2\pi i)=0$.  
This implies $p_2=p_1-p_0$.  Similarly the other edges lead to 
other relations among the $p_i$ and $q_i$, namely:

\begin{tabular}
{l r @{\,\,} c @{\,\,} l @{\qquad\qquad} l  r @{\,\,} c @{\,\,} l}
$z_0z_1$:&$ p_2-p_3+p_4$&$=$&$0$&$z_0z_2$:&$ -p_1+p_3+q_3-p_4-q_4$&$=$&$0$\\
$z_1z_2$:&$ p_0-q_3+q_4$&$=$&$0$&$z_1z_3$:&$ -p_0-q_0+q_2-p_4-q_4$&$=$&$0$\\
$z_2z_3$:&$ q_0-q_1+p_4$&$=$&$0$&$z_2z_4$:&$ -p_0-q_0+p_1+q_1-p_3$&$=$&$0$\\
$z_3z_4$:&$ p_0-p_1+p_2$&$=$&$0$&$z_3z_0$:&$ p_1+q_1-p_2-q_2+q_4$&$=$&$0$\\
$z_4z_0$:&$ -q_1+q_2-q_3$&$=$&$0$&$z_4z_1$:&$ q_0-p_2-q_2+p_3+q_3$&$=$&$0$.
\end{tabular}

\medskip
Elementary linear algebra verifies that these relations are equivalent
to the equations $p_2=p_1-p_0$, $p_3=p_1-p_0+q_1-q_0$, $q_3=q_2-q_1$,
$p_4=q_1-q_0$, and $q_4 = q_2-q_1-p_0$, as in the lifted five term
relation (\ref{five term}).  The lemma thus follows for
$(x_0,\ldots,x_4)\in \FT^+$. It is then true in general by analytic
continuation.
\end{proof}

\section{Definition of $\themap$}

We can now describe the map $\themap\colon H_3(\PSL(2,\C);\Z)\to
\eprebloch(\C)$.  

We shall first recall a standard chain complex for homology of
$G=\PSL(2,\C)$, the chain complex of ``homogeneous simplices for
$G$.''  We will, however, diverge from the standard by using only
non-degenerate simplices, i.e., simplices with distinct vertices ---
we may do this because $G$ is infinite.  

Let $C_n(G)$ denote the free
$\Z$-module on all ordered $(n+1)$-tuples $\langle
g_0,\dots,g_n\rangle$ of distinct elements of $G$.  Define
$\delta\colon C_n\to C_{n-1}$ by \begin{equation*}\delta\langle
g_0,\dots,g_n\rangle=\sum_{i=0}^n(-1)^i\langle g_0,\dots,\hat
g_i,\dots,g_n\rangle.
\end{equation*} Then each $C_n$ is a free $\Z G$-module under
left-multiplication by $G$.  Since $G$ is infinite the sequence
\begin{equation*}\cdots \to C_2\to C_1\to C_0\to Z\to 0\end{equation*}
is exact, so it is a $\Z G$-free resolution of $\Z$.  Thus the chain
complex
\begin{equation*}\cdots\to C_2\otimes_{\Z G}\Z\to C_1\otimes_{\Z G}\Z\to
C_0\otimes_{\Z G}\Z\to 0\end{equation*} 
computes the homology of $G$. Note that $C_n\otimes_{\Z G}\Z$ is the
free $\Z$-module on symbols $\langle g_0\co \ldots\co g_n\rangle$, where the
$g_i$ are distinct elements of $G$ and $\langle
g_0\co \ldots\co g_n\rangle=\langle g_0'\co \ldots\co g_n'\rangle$ if and only if
there is a $g\in G$ with $gg_i=g_i'$ for $i=0,\ldots,n$.

Thus an element of $\alpha\in H_3(G;\Z)$ is represented by a sum
\begin{equation}
\sum
\epsilon_i\langle g_0^{(i)}\co \ldots\co g_3^{(i)}\rangle
\label{homology element}
\end{equation}
of homogeneous $3$-simplices for $G$ and their negatives (here each
$\epsilon_i$ is $\pm1$).  The fact that this is a cycle means that the
2-faces of these homogeneous simplices cancel in pairs.  We choose
some specific way of pairing cancelling faces and form a geometric
quasi-simplicial complex $K$ by taking a 3-simplex $\Delta_i$ for each
homogeneous 3-simplex of the above sum and gluing together 2-faces of
these $\Delta_i$ that correspond to 2-faces of the homogeneous
simplices that have been paired with each other.

We call a closed path $\gamma$
in $K$ a \emph{normal} path if it meets no $0$- or $1$-simplices of $K$
and crosses all $2$-faces that it meets transversally.  When such a
path passes through a $3$-simplex $\Delta_i$, entering and departing
at different faces, there is a unique edge $E$ of the 3-simplex between
these faces. We say the path \emph{passes} this edge $E$.

Consider a choice of combinatorial flattening $\bfw_i$ for each
simplex $\Delta_i$. Then for each edge $E$ of a simplex $\Delta_i$ of
$K$ we have a log-parameter $l_E=l_E(\Delta_i,\bfw_i)$ assigned.
Recall that this log-parameter has the form $\log z + s\pi i$ where
$z$ is the cross-ratio parameter associated to the edge $E$ of simplex
$\Delta_i$ and $s$ is some integer.  We call ($s$ mod $2$) the
\emph{parity parameter} at the edge $E$ of $\Delta_i$ and denote it
$\delta_E=\delta_E(\Delta_i,\bfw_i)$.
\begin{definition}\label{parity}
Suppose $\gamma$ is a normal path in $K$. The \emph{parity} of $\gamma$
is the sum ($\sum_E\delta_E$ modulo $2$) of the parity parameters of
all the edges $E$ that $\gamma$ passes. Moreover, if $\gamma$ runs in a
neighbourhood of some fixed vertex $V$ of $K$, then the {\em
log-parameter for the path} is the sum $\sum_E\pm
\epsilon_{i(E)}l_E$, summed over all edges $E$
that $\gamma$ passes, where:
\begin{itemize}
\item
$i(E)$ is the index $i$ of
the simplex $\Delta_i$ that the edge $E$ belongs to and
$\epsilon_{i(E)}$ is the corresponding coefficient $\pm1$ from
equation (\ref{homology element});
\item
the extra sign $\pm$ is $+$ or $-$ according
as the edge $E$ is passed in a counterclockwise or clockwise
fashion as viewed from the vertex.
\end{itemize}
\end{definition}
\begin{theorem}\label{beta exists}
Choose $z\in \partial\overline\H^3$ such that $g_0^{(i)}z, g_1^{(i)}z,
g_2^{(i)}z, g_3^{(i)}z$ are distinct points for each $i$.  This
defines an ideal hyperbolic simplex shape for each simplex $\Delta_i$
of $K$ and an associated cross ratio $x_i=[g_0^{(i)}z\co  g_1^{(i)}z\co 
g_2^{(i)}z\co  g_3^{(i)}z]$.  There is a way of assigning combinatorial
flattenings $\bfw_i=\ell(x_i;p_i,q_i)$ to the simplices of $K$ such that the
parity of any normal path in $K$ is zero and the log-parameter of any
normal path in any vertex neighbourhood of $K$ is zero.

For such an assignment the element
$\sum_i\epsilon_{i}[x_i,p_i,q_i]\in\eprebloch(\C)$ is independent of choices
and only depends on the original homology class $\alpha$.  We denote
it $\themap(\alpha)$.  Moreover, $\themap(\alpha)\in\ebloch(\C)$ and
$\themap\colon H_3(\PSL(2,\C);\Z)\to\ebloch(\C)$ is a homomorphism.
\end{theorem}

To prove this theorem we will need a general relation (Lemma
\ref{cycle} below) in $\eprebloch(\C)$ 
that follows from the lifted five term relation.

\section{Consequences of the lifted five term relation}\label{consequences}

Let $K$ be a simplicial complex obtained by gluing 3-simplices
$\Delta_1,\dots,\Delta_n$ together in sequence around a common edge
$E$.  Thus, for each index $j$ modulo $n$, $\Delta_j$ is glued
to each of $\Delta_{j-1}$ and $\Delta_{j+1}$ along one of the two
faces of $\Delta_j$ incident to $E$.  Suppose, moreover, that the
vertices of each $\Delta_j$ are ordered such that orderings agree on
the common 2-faces of adjacent 3-simplices.

There is then a sequence $\epsilon_1=\pm1$, $\ldots$,
$\epsilon_n=\pm1$ such that the 2-faces used for gluing all cancel in
the boundary of the 3-chain $\sum_{j=1}^n\epsilon_j\Delta_j$. (Proof:
choose $\epsilon_1=1$ and then for $i=2,\dots,n$ choose $\epsilon_i$
so the common face of $\Delta_{i-1}$ and $\Delta_i$ cancels.  The
common face of $\Delta_n$ and $\Delta_1$ must then cancel since
otherwise that face occurs with coefficient $\pm2$ in
$\partial\sum_{j=1}^n\epsilon_j\Delta_j$, and $E$ occurs with
coefficient $\pm2$ in $\partial\partial\sum_{j=1}^n\epsilon_j\Delta_j$.)

Suppose now further that a combinatorial flattening $\bfw_i$ has been
chosen for each $\Delta_j$ such that the ``signed sum'' of log
parameters around the edge E vanishes and the same for parity
parameters:
\begin{equation}\sum_{j=1}^n\epsilon_j l_E(\Delta_j,\bfw_j)=0,\quad
\sum_{j=1}^n\epsilon_j \delta_E(\Delta_j,\bfw_j)=0 .
\label{around E}\end{equation}

We think of the edge $E$ as being vertical, so that we can label the
two edges other than $E$ of the common triangle of $\Delta_j$ and
$\Delta_{j+1}$ as $T_j$ and $B_j$ (for ``top'' and ``bottom'').  Let
$\bfw_j'$ be the flattening obtained from $\bfw_j$ by adding
$\epsilon_j\pi i$ to the log parameter at $T_j$ and its opposite edge
and subtracting $\epsilon_j\pi i$ from the log parameter at $B_j$ and
its opposite edge.  If we do this for each $j$ then the total log
parameter and parity parameter at any edge of the complex $K$ is not
changed (we sum log-parameters with the appropriate sign
$\epsilon_j$): --- at $E$ no log-parameter has changed while at every
other edge $\pi i$ has been added at one of the two simplices at the
edge and subtracted at the other.

\begin{lemma}\label{cycle}
With the above notation,
\begin{equation*}\sum_{j=1}^n\epsilon_j[\bfw_j]=
\sum_{j=1}^n\epsilon_j[\bfw_j']~~\in\eprebloch(\C),\end{equation*}
where we are using $[\bfw]$ as a shorthand for $[\ell^{-1}\bfw]$ 
(i.e., $[\bfw]$ means $[z,p,q]$ where $\ell(z;p,q)=\bfw$; see
Lemma \ref{Cover is flattenings}).\end{lemma}

\begin{proof}
Each of the simplices $\Delta_i$ has an associated ideal hyperbolic
structure compatible with the combinatorial flattenings $\bfw_j$. This
ideal hyperbolic structure is also compatible with the flattening
$\bfw_j'$.  Choose a realization of $\Delta_1$ as an ideal simplex in
$\overline\H^3$.  We think of this as a mapping of $\Delta_1$ to
$\overline\H^3$.  We can extend this to a mapping of $K$ to $\overline
\H^3$ which maps each $\Delta_j$ to an ideal simplex with shape
appropriate to its combinatorial flattening.  Adjacent simplices will
map to the same side of their common face in $\overline\H^3$ if their
orientations or the signs $\epsilon_j$ do not match and will be on
opposite sides otherwise.  The fact that the signed sums of log and
parity parameters at edge $E$ are zero guarantees that the
identifications match up as we go once around the edge $E$ of $K$.

Note that $K$ has $n+2$ vertices.  We first consider the special case
that $n=3$ and there is an ordering $v_0,\dots,v_4$ of the five
vertices of $K$ that restricts to the given vertex ordering for each
simplex.  We also assume the five vertices of $K$ map to distinct
points $z_0,\dots,z_4$ of $\partial \H^3$.

Each simplex $\Delta_j$ for $j=1,2,3$ has vertices obtained by
omitting one the five vertices $v_0,\dots,v_4$.  Denote by $\Delta_4$
and $\Delta_5$ the simplices obtained by omitting each of the other
two vertices. The fact that the common 2-faces of the $\Delta_j$
cancel when taking boundary of the chain
$\epsilon_1\Delta_1+\epsilon_2\Delta_2+\epsilon_3\Delta_3$ means that,
up to sign this sum corresponds to three summands of the chain
$\partial\langle v_0,\dots,v_4\rangle =\sum(-1)^i\langle
v_o,\dots,\hat{v_i},\dots,v_4\rangle$. Choose $\epsilon_4$ and
$\epsilon_5$ so that $\sum_{j=1}^5\epsilon_j\Delta_j$ is
$\pm\partial\langle v_0,\dots,v_4\rangle$.  

We now claim that we can choose unique combinatorial flattenings
$\bfw_4$ and $\bfw_5$ of $\Delta_4$ and $\Delta_5$ so that the signed
sum of log parameters and parity parameters at any edge of
$K\cup\Delta_4\cup\Delta_5$ is zero.  Indeed, this claim does not
depend on the order of the vertices, so by permuting the vertices we
can assume the five vertices are ordered so that $\Delta_1$,
$\Delta_2$ and $\Delta_3$ are the first three simplices occuring in
the five term relation. Then the common edge $E$ is $v_3v_4$ and the
fact that the simplices fit together around this edge is the condition
that their cross-ratio parameters $x_0$, $x_1$, and $x_2$ satisfy
$x_0x_1^{-1}x_2=1$.  Writing the flattenings as elements of $\Cover$
as $(x_0;p_0,q_0)$, $(x_1;p_1,q_1)$, and $(x_2;p_2,q_2)$, the equation
saying signed sum of log parameters at this edge is zero is $(\log x_0
+ p_0\pi i) -(\log x_1 +p_1\pi i)+(\log x_2 +p_2\pi i)=0$. If $x_0$,
$x_1$, and $x_2$ are in the complex upper half plane this implies the
equation $p_2=p_1-p_0$ of the lifted five term relation, while
otherwise it implies the appropriate analytic continuations in
$\Cover$ of this.  The desired choice of flattenings of $\Delta_4$ and
$\Delta_5$ is thus determined as in the lifted five term relation by
the choice of $p_0$, $p_1$, $q_0$, $q_1$, and $q_2$ (namely
$p_3=p_1-p_0+q_1-q_0$, $q_3=q_2-q_1$, $p_4=q_1-q_0$, and
$q_4=q_2-q_1-p_0$ if $x_0,x_1,x_2$ are in the upper half plane and
otherwise the appropriate analytic continuation).

Note that $\bfw_4$ and $\bfw_5$ do not change if we replace
$\bfw_1,\dots,\bfw_3$ by $\bfw_1',\dots,\bfw_3'$ (with the above
reordering of vertices this just subtracts $1$ from each of $q_0$,
$q_1$ and $q_2$ in the lifted five term relation, so it does not alter
$p_3,q_3,p_4,q_4$).  By Lemma \ref{geometric five term} we then have
\begin{equation*}
\begin{aligned}
\epsilon_1\bfw_1+\epsilon_2\bfw_2+\epsilon_3\bfw_3 &=
-(\epsilon_4\bfw_4+\epsilon_5\bfw_5)\\
\epsilon_1\bfw_1'+\epsilon_2\bfw_2'+\epsilon_3\bfw_3' &=
-(\epsilon_4\bfw_4+\epsilon_5\bfw_5),
\end{aligned}
\end{equation*}
proving this case.

We next consider the case that for some index $j$ modulo $n$ the
images of $\Delta_j$ and $\Delta_{j+1}$ in $\overline\H^3$ do not
coincide, so their union has five distinct vertices.  By cycling our
indices we may assume $j=1$. Since the orderings of the vertices of
$\Delta_1$ and of $\Delta_{2}$ agree on the three vertices they have
in common, there is an ordering of all five vertices compatible with
both $\Delta_1$ and $\Delta_{2}$. Let $\Delta_0$ be the simplex
determined by the common edge $E$ and the two vertices that $\Delta_1$
and $\Delta_{2}$ do not have in common.  Then there is an
$\epsilon_0=\pm1$ such that the common faces of $\Delta_0$, $\Delta_1$,
and $\Delta_{2}$ cancel in the boundary of the chain
$\epsilon_0\Delta_0+\epsilon_1\Delta_1+\epsilon_{2}\Delta_{2}$.
Choose a flattening $\bfw_0$ of $\Delta_0$ such that
$\epsilon_0l_E(\Delta_0,\bfw_0)+\epsilon_1l_E(\Delta_1,\bfw_1)+
\epsilon_{2}l_E(\Delta_{2},\bfw_{2}) =0$.  Then the relation of the
lemma has already been proved for $\bfw_0$, $\bfw_1$, $\bfw_2$, and
by subtracting this relation from the relation to be proved for
$\bfw_1,\dots,\bfw_n$ we obtain a case of the lemma with one less
simplices.  Thus, if we assume the lemma proved for $n-1$ simplices
then this case is also proved.

The above induction argument fails only for the case that there
are $2m$ simplices that alternately ``fold back on each other'' so
that their images in $\overline\H^3$ all have the same four
vertices. The above induction eventually reduces us to this case
(usually with $m=1$).  We must therefore deal with this situation to
complete the proof.  We first consider the case that $m=1$ so $n=2$.
We then have four vertices $z_0,\ldots,z_3$ in
$\partial\overline\H^3$.  We assume the edge $E$ is $z_0z_1$.  Then
the ordering of the vertices of the faces $z_0z_1z_2$ and $z_0z_1z_3$
is the same in each of $\Delta_1$ and $\Delta_2$.  Choose a new point
$z_4$ in $\partial\overline\H^3$ distinct from $z_0,\dots,z_3$ and
consider the ordered simplex with vertices $z_0,z_1,z_2$ ordered as
above followed by $z_4$.  Call this $\Delta_3$.  Similarly make
$\Delta_4$ using $z_0,z_1,z_3$ ordered as above followed by $z_4$.
Choose flattenings of $\Delta_3$ and $\Delta_4$ so that the signed sum
of log parameters for $\Delta_1,\Delta_3,\Delta_4$ around $E$ is
zero.  Then we obtain a three simplex relation of the type already
proved for $\Delta_1,\Delta_3,\Delta_4$ and another for
$\Delta_2,\Delta_3,\Delta_4$, and the difference of these two
relations gives the desired two-simplex relation.

More generally, if we are in the above ``folded'' case with $m>1$ we
can use an instance of the three-simplex relation to replace one of
the $2m$ simplices by two. We then use the induction step to replace
one of these new simplices together with an adjacent old simplex by
one simplex and then repeat for the other new simplex.  In this way we
reduce to a relation involving $2m-1$ simplices, completing the proof.
\end{proof}

Before we continue, we note a consequence of the first case we
considered in the above proof that we will need later.  If vertices
have been reordered as in that proof then the relation we proved can
be written (with the appropriate relationship among $p_0,p_1,p_2$):
\begin{equation}\begin{aligned}{}
&[x,p_0,q_0]-[y,p_1,q_1]+[y/x,p_2,q_2]={}\\
&[x,p_0,q_0-1]-[y,p_1,q_1-1]+
[y/x,p_2,q_2-1].\label{homo}\end{aligned}\end{equation} 
This is true for any choice
of $q_0,q_1,q_2$ so long as $p_0,p_1,p_2$ satisfy the appropriate
relation.  Thus if we just change $q_0$ and subtract the resulting
equation from the above we get \begin{equation*}
[x,p_0,q_0]-[x,p_0,q_0']=[x,p_0,q_0-1]-[x,p_0,q_0'-1].  \end{equation*} From the
versions of the above three-simplex case with different orderings of
the vertices we can derive three versions of this relation:
\begin{equation}
\begin{aligned}{}
[x,p,q]-[x,p,q']&=[x,p,q-1]-[x,p,q'-1]\\
[x,p,q]-[x,p',q]&=[x,p-1,q]-[x,p'-1,q]\label{three equns}\\
[x,p,q]-[x,p+s,q-s]&=[x,p+1,q-1]-[x,p+s+1,q-s-1]
\end{aligned}
\end{equation}
From these
we obtain:
\begin{lemma}\label{super transfer}
$ [x,p,q]=
pq[x,1,1]-(pq-p)[x,1,0] -(pq-q)[x,0,1]+(pq-p-q+1)[x,0,0]$.
\end{lemma}
\begin{proof}
The first of the relations (\ref{three equns}) implies
$[x,p,q]=[x,p,q-1]+[x,p,1]-[x,p,0]$ and applying this repeatedly shows
\begin{equation}[x,p,q]=q[x,p,1]-(q-1)[x,p,0].\label{dd}\end{equation}
The second equation of (\ref{three
equns}) implies similarly that $[x,p,q]=p[x,1,q]-(p-1)[x,0,q]$ and
using this to expand each of the terms on the right of (\ref{dd})
gives the desired equation.
\end{proof}

Up to this point we have not used the transfer relation
(\ref{transfer}). We digress briefly to show that the transfer
relation almost follows from the five term relation.
\begin{proposition}\label{kappa}
If $\eprebloch'(\C)$ and $\ebloch'(\C)$ are defined like
$\eprebloch(\C)$ and $\ebloch(\C)$ but without the transfer relation,
then in $\eprebloch'(\C)$ the element
$\kappa:=[x,1,1]+[x,0,0]-[x,1,0]-[x,0,1]$ is independent of $x$ and
has order $2$. Moreover, $\eprebloch'(\C)={\eprebloch}(\C)\times C_2$
and $\ebloch'(\C)={\ebloch}(\C)\times C_2$, where $C_2$ is the cyclic
group of order 2 generated by $\kappa$.
\end{proposition}
\begin{proof}
If we subtract equation (\ref{homo}) with $p_0=p_1=q_0=q_1=q_2=1$ from
the same equation with $p_0=p_1=0$, $q_0=q_1=q_2=1$ we obtain
$[x,1,1]-[y,1,1]-[x,0,1]+[y,0,1]=[x,1,0]-[y,1,0]-[x,0,0]+[y,0,0]$,
which rearranges to show that $\kappa$ is independent of $x$.  The
last of the equations (\ref{three equns}) with $p=q=0$ and $s=-1$
gives $2[x,0,0]=[x,1,-1]+[x,-1,1]$ and expanding the right side of
this using (ii) gives $2[x,0,0]=-2[x,1,1]+2[x,1,0]+2[x,0,1]$, showing
that $\kappa$ has order dividing 2.

To show $\kappa$ has order exactly 2 we note that there is a
homomorphism $\epsilon:\eprebloch'(\C)\to \Z/2$ defined on generators
by $[z,p,q]\mapsto (pq\text{ mod }2)$.  Indeed, it is easy to check
that this vanishes on the lifted five-term relation, and is thus well
defined on $\eprebloch'(\C)$.  Since $\epsilon(\kappa)=1$ we see
$\kappa$ is non-trivial.
Finally, Lemma \ref{super transfer} implies that the effect of the
transfer relation is simply to kill the element $\kappa$, so the final
sentence of the proposition follows.
\end{proof}

\begin{lemma}\label{1-x}
For any $[x,p,q]\in\eprebloch(\C)$ one has $[x,p,q]+[1-x,-q,-p]=2[1/2,0,0]$.
\end{lemma}
\begin{proof}
Assume first that $0<y<x<1$. Then, as remarked in the proof of
Proposition \ref{R}, \begin{equation*}\begin{aligned}{}
[x,p_0,q_0]&-[y,p_1,q_1]+[\frac
yx,p_1-p_0,q_2]-[\frac{1-x^{-1}}{1-y^{-1}},p_1-p_0+q_1-q_0,q_2-q_1]+{}\\
&[\frac{1-x}{1-y},q_1-q_0,q_2-q_1-p_0]=0\end{aligned}\end{equation*}
is an instance of the lifted five term relation.  Replacing $y$ by
$1-x$, $x$ by $1-y$, $p_0$ by $-q_1$, $q_1$ by $-q_0$, $q_0$ by
$-p_1$, $q_1$ by $-p_0$, and $q_2$ by $q_2-q_1-p_0$ replaces this
relation by exactly the same relation except that the first two terms
are replaced by $[1-y,-q_1,-p_1]-[1-x,-q_0,-p_0]$. Thus subtracting
the two relations gives:
\begin{equation*}[x,p_0,q_0]-[y,p_1,q_1]-[1-y,-q_1,-p_1]+[1-x,-q_0,-p_0]=0.
\end{equation*}
Putting $[y,p_1,q_1]=[1/2,0,0]$ now proves the lemma for $1/2<x<1$.
But since we have shown this as a consequence of the lifted five term
relation, we can analytically continue it over the whole of $\Cover$.
\end{proof}
\begin{proposition}
The following sequence is exact:
\begin{equation*}0\to\C^*\stackrel\chi\longrightarrow
{\eprebloch}(\C)\to\prebloch(\C)\to 0\end{equation*} where
$\eprebloch(\C)\to\prebloch(\C)$ is the natural map and
$\chi(z):= [z,0,1]-[z,0,0]$ for $z\in \C^*$. 
\end{proposition}
\begin{proof}
Denote $\{z,p\}:=[z,p,q]-[z,p,q-1]$ which is independent of $q$ by the
first equation of (\ref{three equns}).  By Lemma \ref{1-x} we have
$[z,p,q]-[z,p-1,q]=-\{1-z,-q\}$.  It follows that elements of the form
$\{z,p\}$ generate $\Ker\bigl(\eprebloch(\C)\to\prebloch(\C)\bigr)$. Computing
$\{z,p\}$ using Lemma \ref{super transfer} and the transfer relation,
one finds $\{z,p\}=\{z,0\}$ which only depends on $z$.  Thus the
elements $\{z,0\}=\chi(z)$ generate
$\Ker\bigl({\eprebloch}(\C)\to\prebloch(\C)\bigr)$.  If we take equation
(\ref{homo}) with even $p_i$ and subtract the same equation with the
$q_i$ reduced by $1$ we get an equation that says that
$\chi\colon\C^*\to \Ker\bigl({\eprebloch}(\C)\to\prebloch(\C)\bigr)$ is a
homomorphism.  We have just shown that it is surjective, and it is
injective because $R\circ\chi$ is the map $\C^*\to \C/\pi^2$ defined
by $z\mapsto \frac{\pi i}2\log z$.
\end{proof}

We can now describe the relationship of our extended groups with the
``classical'' ones.
\begin{theorem} There is a commutative diagram with exact rows and
columns
\begin{equation*}
\begin{CD}
&& 0 &&  0 &&  0 \\
&& @VVV @VVV @VVV \\
0 @>>> \mu^* @>>> \C^* @>>> \C^*/\mu^* @>>> 0\\
&& @V\chi|\mu^* VV @V\chi VV @V\beta VV @VVV\\
0 @>>> \ebloch(\C) @>>> \eprebloch(\C) @>\nu>> \C\wedge\C @>>> 
K_2(\C) @>>> 0\\
&& @VVV @VVV @V\epsilon VV @V=VV\\
0 @>>> \bloch(\C) @>>> \prebloch(\C) @>\nu'>> \C^*\wedge\C^* @>>> 
K_2(\C) @>>> 0\\
&& @VVV @VVV @VVV @VVV\\
&& 0 &&  0 &&  0 && 0
\end{CD}
\end{equation*}
Here $\mu^*$ is the group of roots of unity and the labelled maps
defined as follows:
\begin{equation*}\begin{aligned}
\chi(z)&=[z,0,1]-[z,0,0]\in
\eprebloch(\C);\\
\nu[z,p,q]&=(\log z +
p\pi i)\wedge(-\log(1-z)+ q\pi i);\\
\nu'[z]&=2\bigl(z \wedge (1-z)\bigr);\\
\beta[z]&=\log z \wedge \pi i;\\
\epsilon(w_1\wedge w_2)&= -2(e^{w_1}\wedge e^{w_2});
\end{aligned}\end{equation*}
and the unlabelled maps are the obvious ones.
\end{theorem}
\begin{proof}
The top horizontal sequence is trivially exact while the other two are
exact at their first two non-trivial groups by definition of $\ebloch$
and $\bloch$.  The bottom row is exact also at its other two places by
Milnor's definition of $K_2$.  The exactness of the third vertical
sequence is elementary and the second one has just been proved.  The
commutativity of all but the top left square is elementary. A diagram
chase confirms that $\chi$ maps $\mu^*$ to $\ebloch(\C)$ and that the
left vertical sequence is also exact. Another confirms exactness of
the middle row.
\end{proof}

\section{Proof of Theorem \ref{beta exists}}

We must first recall some notation and results from \cite{neumann}.

The complex $K$ of Theorem \ref{beta exists} is what is called an
``oriented 3-cycle'' in \cite{neumann}.  That is, it is a finite
quasi-simplicial 3-complex such that the complement $K-K^{(0)}$ of the
vertices is an oriented 3-manifold. (``Quasi-simplicial'' means that
$K$ is a CW-complex obtained by gluing together simplices along their
faces in such a way that the interior of each face of each simplex
injects into the resulting complex.)  Each simplex $\Delta_i$ of $K$
has an orientation compatible with the ordering of its vertices, and
this orientation is compatible with or opposite to the orientation of
$K$ according as $\epsilon_i$ is $+1$ or $-1$.

To an oriented $3$-simplex $\Delta$ of $K$ we associate a $2$-dimensional 
bilinear space $J_\Delta$ over $\Z$ as follows.  As a $\Z$-module 
$J_\Delta$ is generated by the six edges $e_0,\dots,e_5$ of $\Delta$ 
(see Fig.~\ref{figdelta}) with the relations: 
\begin{gather*}
e_i - e_{i+3} = 0 \quad\hbox{for }i=0,1,2.\\
e_0 + e_1 + e_2 = 0.
\end{gather*}
\begin{figure}[hbtp]
\centerline{\epsfxsize.3\hsize\epsffile{ebfig3.eps}}
\caption{\label{figdelta}}
\end{figure}

Thus, opposite edges of $\Delta$ represent the same element of 
$J_\Delta$, so $J_\Delta$ has three ``geometric'' generators, 
and the sum of these three generators is zero.  
The bilinear form on $J_\Delta$ is the non-singular skew-symmetric 
form given by
\begin{equation*}
\langle e_0,e_1\rangle = \langle e_1,e_2\rangle = \langle e_2,e_0\rangle = 
-\langle e_1,e_0\rangle = -\langle e_3,e_1\rangle = -\langle e_0,e_2\rangle
= 1. 
\end{equation*}

Let $J$ be the direct sum $\coprod J_\Delta$, summed over the
oriented $3$-simplices of $K$.  For $i=0,1$ let $C_i$ be the free
$\Z$-module on the \emph{unoriented} $i$-simplices of $K$.  Define homomorphisms
\begin{equation*}
\alpha\colon C_0\larrow C_1
\quad\text{and}\quad
\beta\colon C_1\larrow J
\end{equation*}
as follows.
$\alpha$ takes a vertex to the sum of the incident edges 
(with an edge counted twice if both endpoints are at the given vertex).
The $J_\Delta$ component of $\beta$ takes an edge $E$ of $K$ to the 
sum of those edges $e_i$ in the edge set $\{e_0,e_1,\dots,e_5\}$ of $\Delta$ 
which are identified with $E$ in $K$.

The natural basis of $C_i$ gives an identification of $C_i$ with its
dual space and the bilinear form on $J$ gives an identification of $J$
with its dual space.  With respect to these identifications, the dual
map
\begin{equation*}
\alpha^*\colon C_1\larrow C_0
\end{equation*}
is easily seen to map an edge $E$ of $K$ to the sum of its endpoints,
and the dual map 
\begin{equation*}
\beta^*\colon J \larrow C_1
\end{equation*}
can be described as follows.  To each $3$-simplex $\Delta$ of $K$ we 
have a map $j = j_\Delta$ of the edge set $\{e_0,e_1,\dots,e_5\}$ of 
$\Delta$ to the set of edges of $K$: put $j(e_i)$ equal to the edge 
that $e_i$ is identified with in $K$.  For $e_i$ in $J_\Delta$ we have 
\begin{equation*}
\beta^*(e_i) = j(e_{i+1}) - j(e_{i+2}) + j(e_{i+4}) - j(e_{i+5})
\quad\hbox{(indices mod 6).}
\end{equation*}

Let $K_0$ be the result of removing a small open cone neighborhood of 
each vertex $V$ of $K$, so $\partial K_0$ is the disjoint union of 
the links $L_V$ of the vertices of $K$.

\begin{theorem}[\cite{neumann}, Theorem 4.2]
The sequence 
\begin{equation*}
\Jscr\colon\quad0\larrow C_0\larrow^\alpha C_1\larrow^\beta J
\larrow^{\beta^*} C_1\larrow^{\alpha^*} C_0\larrow 0
\end{equation*}
is a chain complex. Its homology groups $H_i(\Jscr)$ (indexing the
non-zero groups of $\Jscr$ from left to right with indices $5,4,3,2,1$) are
\begin{gather*}
H_5(\Jscr)=0,\quad H_4(\Jscr)=\Z/2,\quad H_1(\Jscr)=\Z/2,\\
H_3(\Jscr)=\Hscr\oplus H^1(K;\Z/2),\quad
H_2(\Jscr)=H_1(K;\Z/2),
\end{gather*}
where $\Hscr=\Ker(H_1(\partial K_0;\Z)\to H_1(K_0;\Z/2))$.  Moreover,
the isomorphism $H_2(\Jscr)\to H_1(K;\Z/2)$ results by interpreting an
element of $\Ker(\alpha^*)\subset C_1$ as an unoriented 1-cycle in $K$.
\qed 
\end{theorem}

If $\Delta$ is an ideal simplex and 
$\bfw=(w_0,w_1,w_2)$ is a flattening of it, then denote
\begin{equation*}
\xi(\bfw):=w_1e_0-w_0e_1\in J_\Delta\otimes\C.\end{equation*}
This
definition is only apparantly unsymmetric since
$w_1e_0-w_0e_1=w_2e_1-w_1e_2 =w_0e_2-w_2e_0$. 

Give each simplex $\Delta_i$ of the complex $K$ of Theorem
\ref{beta exists} the flattening
$\bfw_i^{(0)}:=\ell(x_i;0,0)$. Denote by $\omega$ the element of
$J\otimes\C$ whose $\Delta_i$-component is $\xi(\bfw_i^{(0)})$ for
each $i$.  That is, the ${\Delta_i}$-component of $\omega$ is
$-\bigl(\log(1-x_i)e_0+\log(x_i)e_1\bigr)$.  
\begin{lemma}
$\frac1{\pi
i}\beta^*(\omega)$ is an integer class in the kernel of
$\alpha^*$, so it represents an element of the homology group
$H_2(\Jscr)$.  Moreover
 this element in $H_2(\Jscr)$
vanishes, so $\frac1{\pi i}\beta^*(\omega)=\beta^*(x)$ for some
$x\in J$. 
\end{lemma}
\begin{proof}
Let $\overline J_\Delta$ be defined like $J_\Delta$ but without the
relation $e_0+e_1+e_2=0$, so it is generated by the six edges
$e_0,\dots,e_5$ of $\Delta$ with relations $e_i=e_{i+3}$ for
$i=0,1,2$.  Let $\overline J$ be the direct sum $\coprod\overline
J_\Delta$ over 3-simplices $\Delta$ of $K$.

The map $\beta^*\colon J\to C_2$ factors as
\begin{equation*}
\beta^*\colon J\larrow^{\beta_1}\overline J\larrow^{\beta_2}C_2,
\end{equation*}
with $\beta_1$ and $\beta_2$ defined on each component by:
\begin{equation*}
\begin{split}
\beta_1(e_i)&=e_{i+1}-e_{i+2}\\
\beta_2(e_i)&=j(e_i)+j(e_{i+3})
\end{split}
\quad\Biggr\}\text{ for $i=0,1,2$.}
\end{equation*}

Note that $\beta_1(\xi(\bfw))=w_0e_0+w_1e_1+w_2e_2 \in \overline
J_\Delta\otimes\C$.  Thus if $E$ is an edge of $K$ then the
$E$-component of $\beta^*(\omega)=\beta_2\beta_1(\omega)$ is the sum
of the log-parameters for $E$ in the ideal simplices of $K$ around
$E$ and is hence a multiple of $\pi i$.  In fact, it would be a
multiple of $2\pi i$ except for the adjustments by multiples of $\pi
i$ that are involved in forming the log-parameters from the
logarithms of cross-ratio parameters of the simplices.  We claim that
these adjustments add to an even multiple of $\pi i$ around each edge,
so in fact
\begin{equation}\frac1{\pi i}\beta^*(\omega)\in 2C_2,\label{even path}
\end{equation} 
Once this is proved the lemma follows, since the isomorphism
$H_2(\Jscr)\to H_1(K;\Z/2)$ is the map which interprets an element of
$\Ker(\alpha^*)$ as an unoriented 1-cycle in $K$, and equation
(\ref{even path}) says this 1-cycle is zero modulo 2.

In the terminology of Definition \ref{parity} our claim can be
restated that the parity of the normal path that circles the edge $E$
is zero for each edge $E$.  In fact, the parity of any normal path in
$K$ is zero.  To see this, note that as we follow a normal path the
contribution to the parity as we pass an edge of a simplex is $0$ if
the edge is the $01$, $03$, $12$, or $23$ edge of the simplex and the
contribution is $\pm1$ if it is the $02$ or $13$ edge.  Consider the
orientations of the triangular faces we cross as we follow the path,
where the orientation is the one induced by the ordering of its
vertices. As we pass a $02$ or $13$ edge this orientation changes
while for the other edges it does not.  Since $K$ is oriented, we must
have an even number of orientation changes as we traverse the normal
path, proving the claim.  (This argument, which simplifies my original
one, is due to Brian Bowditch.)
\end{proof}

Let $\omega':=\omega-\pi i x\in J\otimes\C$ with $x$ as in the lemma,
so $\beta^*(\omega')=0$.  The $\Delta_i$-component of $\omega'$ is
$\xi(\bfw_i)$, where $\bfw_i=\ell(x_i;p_i,q_i)$ for suitable integers
$p_i,q_i$.  These integers $p_i$ and $q_i$ are the coefficients
ocurring in the element $x\in J$, which is only determined by the
lemma up to elements of $\Ker(\beta^*)$.  We claim that for suitable
choice of $x$, the $\bfw_i$ satisfy the parity and log-parameter
conditions of first paragraph of Theorem \ref{beta exists}.  To see
this we need to review a computation of $H_3(\Jscr)$ from
\cite{neumann}.

We define a map $\gamma'\colon H_3(\Jscr)\to H^1(\partial
K_0;\Z)=\Hom(H_1(\partial K_0),\Z)$ as follows.  Given elements $a\in
H_3(\Jscr)$ and $c\in H_1(\partial K_0)$ we wish to define
$\gamma'(a)(c)$.  It is enough to do this for a class $c$ which is
represented by a normal path $C$ in the link of some vertex of $K$.
Represent $a$ by an element $A\in J$ with
$\beta^*(A)=0$ and consider the element $\beta_1(A)\in \overline J$.
This element has a coefficient for each edge of each simplex of $K$.
To define $\gamma'(a)(c)$ we consider the coefficients of $\beta_1(A)$
corresponding to edges of simplices that $C$ passes and sum these
using the orientation conventions of Definition \ref{parity}.  It is
easy to see that the result only depends on the homology class of $C$.

We can similarly define a map $\gamma_2'\colon H_3(\Jscr)\to
H^1(K_0;\Z/2)=\Hom(H_1(K_0),\Z/2)$ by using normal paths in $K_0$ and
taking modulo 2 sum of coefficients of $\beta_1(A)$.

\begin{lemma}[\cite{neumann}, Theorem 5.1]
The sequence
\begin{equation*}
0\to H_3(\Jscr)\Mapright{30}{(\gamma',\gamma_2')} H^1(\partial
K_0;\Z)\oplus H^1(K_0;\Z/2)\Mapright{25} {r-i^*}H^1(\partial
K_0;\Z/2)\to 0
\end{equation*}
is exact, where $r\colon H^1(\partial K_0;\Z)\to
H^1(\partial K_0;\Z/2)$ is the coefficient map and the map $i^*\colon
H^1(K_0;\Z/2)\to H^1(\partial K_0;\Z/2)$ is induced by the inclusion
$\partial K_0\to K_0$.
\qed\end{lemma}

Returning to the choice of $x$ above, assume we have made a choice so
that the resulting flattenings $\bfw_i$ do not lead to zero
log-parameters and parities for normal paths. Taking $\frac1{\pi i}$
times the log-parameters of normal paths leads as above to an element
$c\in H^1(\partial K_0;\Z)$. Similarly, parities of normal paths leads
to an element of $c_2\in H^1(K_0;\Z/2)$. These elements satisfy
$r(c)=i^*(c_2)$.  The lemma thus gives an element of $H_3(\Jscr)$ that
maps to $(c,c_2)$, and subtracting a representative for this element
from $x$ gives the desired correction of $x$ so the log-parameters and
parities of normal paths with respect to the corresponding changed
$\bfw_i$'s are zero. This completes the proof of the first paragraph
of Theorem \ref{beta exists}.

Suppose now that we have a different choice of flattenings $\bfw_i$
satisfying the parity and log-parameter conditions of the Theorem.
Then the above lemma implies that the difference between the
corresponding elements $x$ represents $0$ in $H_3(\Jscr)$ and is thus
in the image of $\beta$. For $E\in C_2$ the effect of changing $x$ by
$\beta(E)$ is to change the element $\sum_i\epsilon_i [x_i,p_i,q_i]
\in\eprebloch(\C)$ by the corresponding relation of Lemma \ref{cycle}.
Since this is a consequence of the lifted five term relations, the
element in $\eprebloch(\C)$ is unchanged.

Finally we need to show that none of the other choices we made in
defining the element $\lambda(\alpha)$ in Theorem \ref{beta exists}
have an effect.  These were:
\begin{itemize} 
\item the representative of the homology class $\alpha$;
\item the choice of pairing of faces of simplices of $\alpha$ to form the complex
$K$;
\item the choice of the point $z\in\partial\overline\H^3$ in Theorem
\ref{beta exists}.
\end{itemize}
To prove this we will use a bordism theory based on $n$-cycles which
we describe next.

An \emph{ordered $n$-cycle with boundary} is defined by the following:
\begin{itemize}
\item a finite collection of ordered $n$-simplices $\Delta_i$ is given;
\item each simplex has a sign $\epsilon_i=\pm1$;
\item the $(n-1)$ faces of the simplices $\Delta_i$ are given signs as
follows: the sign of the $(n-1)$-face that omits the
$j$-th vertex of $\Delta_i$ is $(-1)^j\epsilon_i$;
\item a subset of the $D$ of the $(n-1)$-faces of the collection of
simplices is given;
\item the $(n-1)$ faces which are not in $D$ are paired in such a way
the the two simplices in each pair have opposite sign.
\end{itemize} 
We can form a geometric $n$-complex $K$ by realizing the $\Delta_i$ by
geometric $n$-simplices and gluing (by affine isomorphisms that
respect the vertex orderings) faces that are paired.  The result is a
space that is an oriented $n$-manifold with boundary in the complement
of its $(n-2)$-skeleton.  We shall use the same letter $K$ to denote
the underlying combinatorial object and the geometric realization.

If $K$ is an ordered $n$-cycle with boundary as above, we say it is
\emph{closed} if $D$ is empty.  If $K$ is not closed, then $D$
inherits the structure of a closed ordered $(n-1)$-cycle, and we
denote this structure by $\partial K$.  The pairing on the faces of
the simplices of $D$ can be described as follows.  Consider the set of
pairs $(\sigma,\tau)$ consisting of an $(n-2)$-face $\tau$ of an
$(n-1)$-face $\sigma$ of a simplex of $K$.  Construct a graph with
this set as vertex set and with edges of two types: $(\sigma,\tau)$ is
connected to $(\sigma',\tau')$ by an edge if either $\sigma$ and
$\sigma'$ are adjacent faces of an $n$-simplex of $K$ and $\tau=\tau'$
is their common $(n-2)$ face, or if $\sigma$ and $\sigma'$ are paired
$(n-1)$-simplices and $\tau$ and $\tau'$ are corresponding faces under
an order-preserving matching of $\sigma$ with $\sigma'$.  Each vertex
$(\sigma,\tau)$ of this graph has valency $2$ except when $\sigma$ is
in $D$, in which case the valency is $1$. Thus the graph consists of a
collection of arcs and cycles.  The endpoints of each arc gives two
$(n-2)$-faces of elements of $D$ that are to be paired.

\section{The true extended Bloch group}

Recall that $\Cover$ consists of four components $X_{00}$, $X_{01}$,
$X_{10}$, and $X_{11}$, of which $X_{00}$ is naturally the universal
abelian cover of $\C-\{0,1\}$.

Let $\liftFT_{00}$ be $\liftFT\cap (X_{00})^5$, so
$\liftFT_{00}=\liftFT_0+2V$ in the notation of Definition
\ref{liftFT}.  As mentioned earlier, $\liftFT_{00}$ is, in fact, equal
to $\liftFT_0$, but we do not need this.

Define $\eeprebloch(\C)$ to be the free $\Z$-module on $X_{00}$
factored by all instances of the relation
\begin{equation}
\sum_{i=0}^4(-1)^i(x_i;2p_i,2q_i)=0\quad\text{with
$\bigl((x_0;2p_0,2q_0),\dots,(x_4;2p_4,2q_4)\bigr)\in\liftFT_{00}$}.
\end{equation}

As before, we have a well-defined map
\begin{equation*}
\nu\colon\eeprebloch(\C)\to\C\wedge\C, 
\end{equation*}
given by $\nu[z;2p,2q]=(\log z+2p\pi i)\wedge(-\log(1-z)+2q\pi i)$,
and we define
\begin{equation}
\eebloch(\C):=\Ker\nu.
\end{equation}

The proof of Proposition \ref{R} shows:
\begin{proposition}
The function 
$R(z;2p,2q):=\calR(z)+\pi i (p\log(1-z)+q\log
z)-\frac{\pi^2}6$
gives a well defined map $X_{00}\to\C/2\pi^2\Z$ and induces a
homomorphism
$R\colon \eeprebloch(\C)\to\C/2\pi^2\Z$.
\qed\end{proposition}

We can repeat the computations in Sect.~\ref{consequences}
word-for-word, replacing anything of the form $[x,p,q]$ by
$[x,2p,2q]$, to show that $\Ker\bigl(\eeprebloch(\C)\to\prebloch(\C)\bigr)$ is
generated by elements of the form $\hat\chi(z):=[z,0,2]-[z,0,0]$.  We
thus get:
\begin{proposition} 
The following sequence is exact:
\begin{equation*}
0\to\C^*\stackrel{\hat\chi}\longrightarrow
\eeprebloch(\C)\to\prebloch(\C)\to 0
\end{equation*} where
$\eprebloch(\C)\to\prebloch(\C)$ is the natural map and
$\hat\chi(z):= [z,0,2]-[z,0,0]$ for $z\in \C^*$. 
\end{proposition}
\begin{proof}
The only thing to prove is the injectivity of $\hat\chi$ which follows
by noting that $R\bigl(\hat\chi(z)\bigr)=\pi i\log z\in\C/2\pi^2$.\end{proof}
\begin{corollary}
We have a commutative diagram with exact rows and columns:
\begin{equation*}
\begin{CD}
&& 0 &&  0  \\
&& @VVV @VVV  \\
&& \Z/2 @>=>> \Z/2 \\
&& @VVV @VVV  \\
0 @>>> \C^* @>\hat\chi>> \eeprebloch(\C) @>>> \prebloch(\C) @>>> 0\\
&& @VV z\mapsto z^2V @VVV @VV=V \\
0 @>>> \C^* @>\chi>> \eprebloch(\C) @>>> \prebloch(\C) @>>> 0\\
&& @VVV @VVV \\
&& 0 &&  0 && 
\end{CD}
\end{equation*}
and analogously for the Bloch group:
\begin{equation*}
\begin{CD}
&& 0 &&  0  \\
&& @VVV @VVV  \\
&& \Z/2 @>=>> \Z/2 \\
&& @VVV @VVV  \\
0 @>>> \mu^* @>\hat\chi>> \eebloch(\C) @>>> \bloch(\C) @>>> 0\\
&& @VV z\mapsto z^2V @VVV @VV=V \\
0 @>>> \mu^* @>\chi>> \ebloch(\C) @>>> \bloch(\C) @>>> 0\\
&& @VVV @VVV \\
&& 0 &&  0 && 
\end{CD}
\end{equation*}
\end{corollary}

\end{document}